\documentclass{article}
\usepackage{authblk}
\usepackage{cite}
\usepackage{amsmath,amssymb,amsfonts}
\usepackage{algorithmic}
\usepackage{graphicx}
\usepackage{textcomp}
\usepackage[frozencache,cachedir=_minted-osmses_2022-mtress_3]{minted}
\usepackage{siunitx}
\usepackage{eurosym}
\usepackage{xcolor}

\usepackage[version=4]{mhchem}

\usepackage{siunitx}
\DeclareSIUnit{\eur}{\text{\euro}}

\begin{document}

\title{MTRESS 3.0 -- Modell Template for Residential Energy Supply Systems
\thanks{This work has been funded by the German Federal Ministry of Education and Research, Grant No. 03SF0624L.}
}

\author[1]{
	Patrik Schönfeldt\thanks{patrik.schoenfeldt@dlr.de}
}
\author[1]{
        Sunke Schlüters\thanks{sunke.schlueters@dlr.de}
}
\author[1]{
        Keno Oltmanns\thanks{keno.oltmanns@dlr.de}
}
\affil[1]{Institute of Networked Energy Systems\\German Aerospace Centre (DLR)\\Oldenburg, Germany}

\maketitle

\begin{abstract}
MTRESS is a tool that facilitates the creation of models
of residential energy supply systems by providing a template with meaningful presets.
This model can then be used to linearly optimise the operation of the energy system.
Version~3.0 enables multiple locations
belonging to one energy system to be defined.
Furthermore, it adds hydrogen as an energy carrier with the possibility to convert
power to gas using electrolysis and to store gas at different pressure levels.
Additionally, we reworked the Python API, giving more flexibility to the user.
\end{abstract}

\section{Introduction}

Graph based formulation of an energy system in order to
formulate a (mixed-integer) linear optimisation problem
can be considered a standard approach.
In Python, there are several tools implementing that approach.
For example, there are
PyPSA~\cite{Brown_2018},
urbs~\cite{urbs_zenodo},
OSeMOSYS~\cite{OSeMOSYS_HOWELLS20115850},
PowerGAMA~\cite{PowerGAMA},
Minpower~\cite{Minpower_6344667},
MOST~\cite{MOST_6670209},
Calliope~\cite{Calliope_Pfenninger2018}, and
oemof.solph~\cite{Krien2020}.
It is also noteworthy that they have similar dependencies,
i.e. most are based on Pyomo~\cite{bynum2021pyomo, hart2011pyomo}.
As a result,
improvements to one of those tools can be adapted for the others.

The present contribution aims to further simplify the creation
of this kind of models, even when the level of detail is increased.
Also, it offers different styles of creating energy system models,
either by preparing formatted data in text files or by scripting in Python.
The implementation of these styles is influenced by the aforementioned
tools Calliope and oemof.solph.
Calliope emphasises a strict separation of data and program, and
considers the definition of a model as data that is parameterised using
text files in YAML and CSV formats.
On the other hand, oemof.solph implements a "create and add" workflow,
so that implementation of the specific model and the data are located
closely to one another in the code.

\section{Architecture of MTRESS 3.0}

The aim of MTRESS is the simple creation of energy supply system models
while keeping model accuracy at a considerably high level.
To give meaningful presets, a number of assumptions is made that are often met.  
With MTRESS~3.0, a number of new concepts is introduced to allow for more flexibility compared to previous versions.
Further -- as MTRESS wraps oemof.solph --
energy systems created using the template can also be extended and further refined using the full flexibility of oemof.solph.

The overall structure is created by building blocks of the five classes
\texttt{Energy\-System},
\texttt{Location}, 
\texttt{EnergyCarrier}, 
\texttt{Component}, and 
\texttt{Demand}.
One major simplification compared to other methods to create
energy system models is that the connections do not have to be
defined by the user.

\subsection{Energy System}

The \texttt{EnergySystem} serves as a container for the model.
It holds global information, e.g. about the time, and defaults
that can be overwritten for specific locations,
e.g. weather data.
Furthermore, all locations of the model and the connections between them
are registered with the \texttt{EnergySystem} object.
Once the energy system is about to be solved,
it makes sure every location has all the necessary connections
and constraints set.

\subsection{Location}

A \texttt{Location} collects different energy carriers, components and demands.

To allow for automatic connections between the components and demands,
every energy carrier (e.g. electricity or heat) and every component
(e.g. a heat pump) can only be defined once per location (or left out).
To define multiple instances of one energy carrier with different
configurations, multiple locations have to be defined.

\subsection{Energy Carrier}
\label{sec: energy carrier}

Each \texttt{EnergyCarrier} represents a form of energy in a energy systems location
and acts as a distribution hub between components.
Currently, electricity, heat, natural gas, and hydrogen are implemented.
These energy carriers may have a ``grid connection'',
meaning that energy may be purchased or sold.
Currently, the following energy carriers are implemented:

\paragraph{Electricity}
Electricity is implemented in a way that allows the tracking of
which energy is used where.
While there is no physical difference, the model distinguishes between electricity
that is produced locally and electricity coming from the grid.
This enables different fees and tariffs based on the origin
of the electricity to be considered.

\paragraph{Heat}
The temperature has a significant impact on the performance
of renewable energy supply systems, e.g. for the Coefficient of Performance (COP) of heat pumps
or the output of a solar thermal collector.
Thus, the energy carrier heat
allows to optimise both temperature and heat.
This is done by defining several discrete temperature levels~\cite{9768967}.
Besides the temperature levels, a reference temperature can be defined.
This can be useful to simplify the model by setting the reference to a
return temperature, resulting in the corresponding return flow to
be considered at zero energy.
Some temperatures, i.e. the ones of sources for heat pumps,
are not considered by the energy carrier.
To emphasise that fact, these sources are defined as anergy sources,
which are not connected to the energy carrier but only to the heat pump.

\paragraph{Gas}
The concept of the energy carrier heat
is expanded to (compressible) gas.
It can be expanded from high pressure to lower
pressure or compressed (using a compressor).
The pressure level plays a role, especially
when energy is stored in gaseous form. 

As hydrogen was  recently implemented into MTRESS in the form of an electrolyzer and a hydrogen compressor, the gas carrier currently finds application in that field.

\subsection{Component}
A \texttt{Component} represents a technology used in the system.
They can serve as sources or connect either energy carriers to each other or different temperature levels (in case of the heat carrier) or pressure levels (in case of a gaseous carrier) of the carrier.
The following collection of components should give an idea of the concept:

\paragraph{Heat Pump}
The heat pump uses electricity and provides heat
(possibly at different temperature levels). Further, it connects to
every available anergy source at the location.
The COPs are automatically calculated based on the information
given by the heat carrier and the anergy sources.

\paragraph{Anergy Sources}
There can be several anergy sources, e.g. down hole heat exchangers
or air source heat exchangers.
They hold a time series for both the temperature and the power limit
that can be drawn from the source.
Furthermore, additionally a total limit can be defined.
This is particularly important for geothermal sources that need to recover.

\paragraph{Combined Heat and Power}
Combined heat and power (CHP) units consume gas and supply both
heat and electricity.
It is possible to choose between different models that linearize
the CHP to consider (or neglect) nonlinear part-load efficiencies.

\paragraph{Renewable Electricity Source}
Renewable electricity sources are representing photovoltaic (PV),
wind power, or river-flow water generators.
They contain a time series for the maximum power that can be drawn
from them.
For PV, the time series can be automatically calculated
by MTRESS from a weather time series.

\paragraph{Energy storage}
Energy storage is implemented per carrier.
Thus, there can be electricity storage (battery),
heat storage, and gas storage.
The three types incorporate the properties of the corresponding energy carrier.
The battery is represented by a state of charge.
We decided to not allow changing voltage levels for the electricity,
as this technical option is typically not implemented.
For heat storage and gas storage,
stored content depends on an optimisation variable
(temperature or pressure),
the used linearization is described in
Sec.~\ref{sec: storage linearization}.

\paragraph{PEM Electrolyzer}
In regard to the generation of hydrogen, a protone-exchange-membrane (PEM) electrolyzer was implemented. It uses electricity and water (\ce{H2O}) to produce hydrogen (\ce{H2}) and oxygen (\ce{O}). For the simplicity of the model, water-input as well as oxygen-output can be neglected. However, the usage of waste heat generated by the electrolyzer was implemented into the model. Usually waste-heat of PEM electrolyzers is  produced at \SI{77}{\celsius}~\cite{Tiktak.2019}, while hydrogen is produced at a pressure of around \SI{30}{\bar}~\cite{htec_datasheet}.

\paragraph{Compressor}
To consider the necessity for different pressure levels of hydrogen, e.g. for transport or usage at the maximum density of \SI{700}{\bar} or a storing at \SI{350}{\bar} for rocket science applications~\cite{h2orizon_2022}, a compressor was implemented. The compressor uses electricity to rise the pressure level of the produced hydrogen. 

\subsection{Demand}

Demands contain time series of energy that is needed.
Depending on the type of demand (e.g. electricity or heat),
it automatically connects to the corresponding energy carrier.
Also, a name identifying the demand has to be given
that is unique for the location.
This is because multiple demands of one type can exist for one location.
The different types have different complexity:
Electricity demand does not need any further specification,
heat and gas demand need a specified temperature or pressure level,
respectively.
Further, energy from electricity and the gaseous carriers is just consumed,
heat demands have a returning energy flow.

\subsection{Helpers}

Besides the part of the energy system, MTRESS offers some helpers
to facilitate working with the models.

\paragraph{Reproducible runs}
It is possible to create an energy system model using the Python API
and to export it to a YAML file afterwards.
This way, models created dynamically using the
Python API can be easily checked and results can be reproduced
without running the full program again.

\paragraph{Analysis of the results}
Energy flows in MTRESS are categorised using tags.
It is possible to filter for these tags, e.g. to sum up
heat coming from heat pumps in several buildings.
Further analysis and plotting can be done by exporting
results to Pandas~\cite{mckinney2010data}.

\paragraph{Energy system graph plotting}
For quick checks and for understanding the layout,
we implemented a simple graph plotting routine.
It displays every node of the final energy system graph
and the connections between those.
Trapezoid shapes signify sources (smaller side up)
or sinks (smaller size down), cycles signify busses
and octagonal shapes signify more complicated energy converters.
Adding to the same functionality in solph,
the MTRESS plotter groups these nodes corresponding to the location
and further according energy carrier, component, or demand they belong.
An example is displayed in Fig.~\ref{fig: example energy system graph}.
For a later publication, the plotted graphs can be exported as \texttt{graphml}
and edited manually in graph visualisation tools.

\section{Linearized Formulation of Energy Storage}
\label{sec: storage linearization}

For the gas storage, available pressure depends
on the storage content. The latter, however,
is an optimisation variable.
The same argument is valid for a (sensible) heat storage when
it comes to the temperature of the medium
that can be taken from the storage.
Thus, linear formulation of such storage types is
not trivial.
As discussed in Sec.\ref{sec: energy carrier},
we implement an approach that uses discrete levels.
Using the example of a gas storage,
it has varying pressure that defines the content of the storage, i.e.
\begin{equation}
    E(t) = c_\mathrm{E} \cdot p(t),
\end{equation}
where $E$ denotes the energy content of the storage at time $t$,
$p$ the pressure inside the storage at time $t$
and $c_E$ is a storage specific constant depending on the energy density of
the gas and the volume of the storage.
Following the definition of the \texttt{EnergyCarrier},
the pressure is discretized
\begin{equation}
    p_n ,\, n \in \{1,\dotsc, N\},
\end{equation}
and \(E_n := c_\mathrm{E} \cdot p_n\) denote the energy contents of the storage
at the defined pressure levels.

Now, we want that an active (or usable) level is signified by the binary
status variable $y_{n}(t) \in \{0, 1\}$ with
\begin{subequations}
    \begin{align}
        y_{n}(t) &= 0 \text{ if } E(t) < E_n
        ,\, \text{and}
        \label{eq: y_n inactive}\\
        y_{n}(t) &=1 \text{ if } E(t) \ge E_n.
        \label{eq: y_n active}
    \end{align}%
    \label{eq: y_n cases}%
\end{subequations}%
Note that \(E\) is being optimised,
so Eq.~\eqref{eq: y_n cases} cannot be read as a causal relation.
We also need a linear formulation.
We suggest
\begin{subequations}
\begin{align}
    y_{n}(t) &\le \frac{E(t)}{E_n},\label{eq: y_n possible}\\
    \hat{y}_{n}(t) &\ge \frac{E(t) - E_n}{E_\mathrm{max}},\\
    \bar{y}_{n}(t) &= 1 - \hat{y}_{n}(t),
    \label{eq: anti y_n}
\end{align}
where \(E_\mathrm{max}\) is the maximum storage content.
Equation~\eqref{eq: y_n possible} guarantees Eq.~\eqref{eq: y_n inactive}
but relaxes Eq.~\eqref{eq: y_n active}
in the sense that \(y_{n}\) is not forced to be active.
To compensate for that,
Eq.~\eqref{eq: anti y_n}, enforces \(\bar{y}_{n} = 0\) for the given case.
The symbol is chosen to emphasise that it can be read as an inverse status.
If this should be strictly the case, also
\begin{equation}
    1 = {y}_{n}(t) + \bar{y}_{n}(t)
\end{equation}
\end{subequations}
has to be defined to eliminate the possibility that
\({y}_{n}(t) = \bar{y}_{n}(t)\), i.e. at \(E = E_n\),
where both can be 1.

Now, note that gas leaving the storage needs to have lower pressure
than the storage and higher pressure is needed to increase the storage content.
Thus,
\begin{subequations}
    \begin{align}
        P_{\mathrm{out}, n}(t)
            &\le y_{n}(t) \cdot P_{\mathrm{out, max}, n},\\
        P_{\mathrm{in}, n}(t)
            &\le \bar{y}_{n}(t) \cdot P_{\mathrm{in, max}, n},
    \end{align}%
    \label{eq: power limit}%
\end{subequations}%
where $P_{\mathrm{out},n}$ and $P_{\mathrm{in},n}$ denote the power flow out of and into the storage, respectively,
at the pressure $p_{n}$. $P_{\mathrm{out},\mathrm{max},n}$ and $P_{\mathrm{in},\mathrm{max},n}$ are the absolute
limits of these flows.
This way, it is guaranteed that \(P_\mathrm{out}(t) = 0\) if the storage content
is not sufficient and
\(P_\mathrm{in}(t) = 0\) if the storage content is too high.

So, at each point in time, the storage content \(E\) defines the pressure \(p\)
in the storage and thus limits the withdrawal pressure or the feed pressure.
As multiple output and input flows can be active at the same time,
it is meaningful to also define a weighted limit
\begin{equation}
    \sum_n P_{\mathrm{out}, n}(t) \cdot w(p_n) \le P_{\mathrm{out}, \mathrm{max}},
\end{equation}
so that the total power flow out of the storage, and analogously for the input flows, is also constrained.

\begin{figure}[bt]
\centerline{\includegraphics[width=\linewidth]{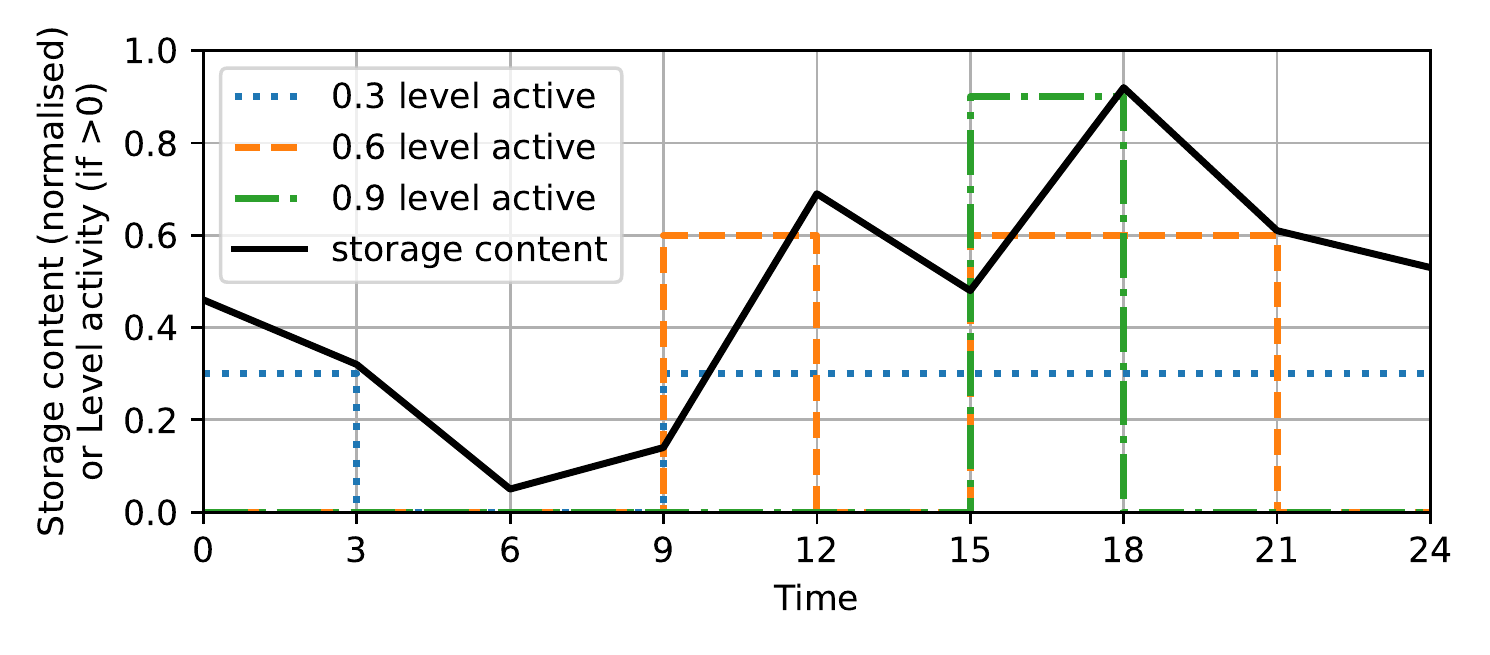}}
\caption{Working principle of the status variable \(y_n\):
If the stored energy exceeds the limit (at the end of a time step),
the corresponding level might be active.}
\label{fig: status example}
\end{figure}

A further complication comes up due to discrete time steps.
A naive implementation might allow to pass a pressure level \(p_n\)
within one time step using \(P_{\mathrm{out}, n}(t) > 0\) -- 
even a full storage can be completely emptied
using the highest pressure level.
If either power limits as in Eq.~\eqref{eq: power limit} are set low
or the time resolution is high, this will not impose a problem.
However, there is a solution that explicitly solves the issue.
\begin{subequations}
    \begin{align}
        P_{\mathrm{out}, n,t}
            &\le y_{n,t} \cdot P_{\mathrm{out, max}, n},\\
        P_{\mathrm{in}, n,t}
            &\le \bar{y}_{n,t} \cdot P_{\mathrm{in, max}, n},
    \end{align}%
    \label{eq: power limit, discrete time}%
\end{subequations}
is a time-discrete reformulation of Eq.~\eqref{eq: power limit},
there the index \(t\) denotes the time interval between \(t\) and \(t+\Delta t\).
If now the constraint for the status variable is defined by the energy
content at the end of that time interval
\begin{equation}
    y_{n,t} \le \frac{E_{t+\Delta t}}{E_n},
\end{equation}
it is no longer possible to cross a level \(p_n\) using the
corresponding power flows.
For a storage with levels at \(E_n = [0, 0.3, 0.6, 0.9, 1] \cdot E_\mathrm{max}\),
the principle is depicted in Fig.~\ref{fig: status example}.
Lines for the highest and the lowest level are omitted because
they are always active or inactive, respectively.
This implies that at any time the lowest level cannot be used
for storing energy and energy at the highest level cannot be obtained
from the storage.
For a fully mixed heat storage, this implication can be read as the fact that
a storage will never reach (exactly) the temperature that is used to feed it.

Generally, however, inside a heat storage stratification is possible.
Thus, a second heat storage model is offered~\cite{9768967},
that works without binary variables.
It models different temperature layers of variable height
that share the same volume.
We assume that the warmest temperatures are at the top and temperature
decreases until the bottom of the storage tank.
This way, multiple temperatures can be served at the same time
using the same storage.

For electricity, the presented method can be used to implement
charging-rates that depend on the state of charge.

\section{Example}

\begin{figure*}[htbp]
\centerline{\includegraphics[width=\linewidth]{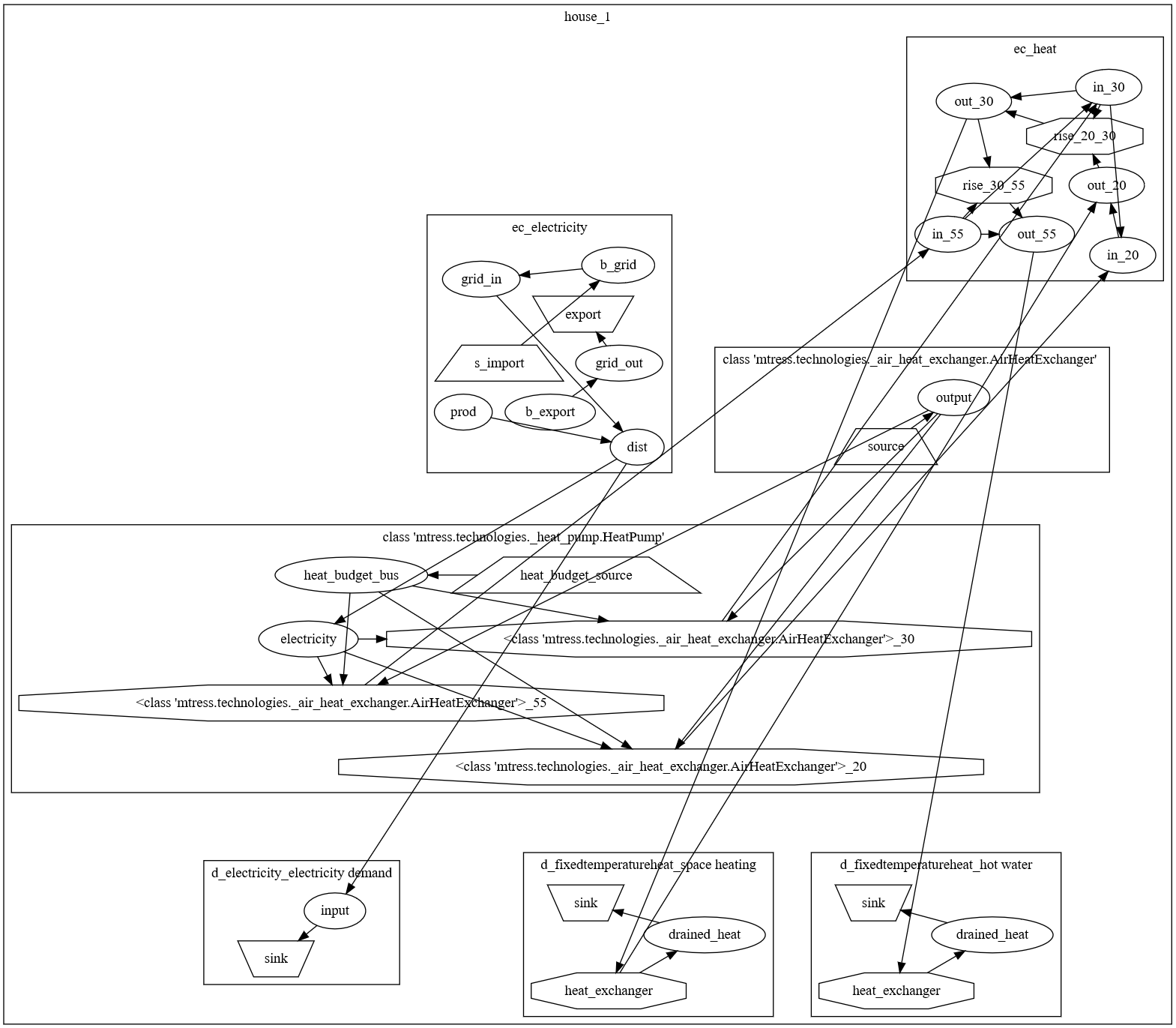}}
\caption{Example energy system graph as plotted by the automatic plotting routine.}
\label{fig: example energy system graph}
\end{figure*}

The case of a single-family home using a heat pump for heat supply will serve
as a basic example.
First of all, the \texttt{EnergySystem} object has to be created.
As stated above, it will nest the locations and contains global
information about the problem.
The simplest way to define the time horizon of the optimisation
is by handing the needed information as a dictionary:
\begin{minted}[xleftmargin=\parindent,breaklines,linenos,firstnumber=1]{Python}
energy_system = EnergySystem(time_index={
    "start": "2021-07-10 06:00:00",
    "end": "2021-07-10 08:00:00",
    "freq": "60T",
})
\end{minted}

Note that the boundaries given above define two time intervals. 
Next, we create an empty location.
As there can be many locations, the location has to be given
a unique name as an identifier.
\begin{minted}[xleftmargin=\parindent,breaklines,linenos,firstnumber=6]{Python}
house = Location(name="SFH")
energy_system.add(house)
\end{minted}

We add the location to the model at this point,
just to have it done.
Generally speaking,
the order of the individual steps of creating the model is not important for the result.
In particular, it is also possible to create and populate the locations
first, before creating the \texttt{EnergySystem} object.

\begin{minted}[xleftmargin=\parindent,breaklines,linenos,firstnumber=8]{Python}
house.add_carrier(
    carriers.Electricity(costs={
        "working_price": 35,  # ct/kWh
        "demand_rate": 0
    })
)
house.add_demand(
    demands.Electricity(
        name="electricity demand",
        time_series=[7, 8.4],  # kW
    )
)
\end{minted}

By giving a working price,
it is implied that energy for that carrier can be purchased
(from an external grid).
The convention for the units is not binding,
e.g. using \SI{}{\eur\per\mega\watt} and \SI{}{\mega\watt}
would also be possible, 
as long as it is consistently done for all carriers and demands.
As for all demands, the electricity demand time series
is given in units of power.
The energy system defined up to this point describes a complete
electricity-only energy system and could be used as such.

For the heat sector, temperature levels and a reference temperature
have to be defined.
The reference temperature has to be lower than or equal to the lowest
return temperature.
Besides, every temperature required by the demands has to be defined
as a level for the carrier.
\begin{minted}[xleftmargin=\parindent,breaklines,linenos,firstnumber=20]{Python}
house.add_carrier(
    carriers.Heat(
        temperature_levels=[
            20, 30, 55],  # °C
        reference_temperature=10,  # °C
    )
)
house.add_demand(
    demands.FixedTemperatureHeat(
        name="space heating",
        flow_temperature=30,  # °C
        return_temperature=20,  # °C
        time_series=[13.37, 42],  # kW
    )
)
house.add_demand(
    demands.FixedTemperatureHeat(
        name="hot water",
        flow_temperature=55,  # °C
        return_temperature=10,  # °C
        time_series=[0, 12],  # kW
    )
)
\end{minted}

Note that we did not define costs for the heat carrier,
an equivalent expression would be \texttt{costs=None}.
Defining it would be interpreted as a connection to the heating grid,
i.e. setting it to zero would make it free.
So, we need an energy conversion technology.
In our example, an air-source heat-pump is used.
It consists of two individual components,
the actual \texttt{HeatPump} and the corresponding
anergy source:
\begin{minted}[xleftmargin=\parindent,breaklines,linenos,firstnumber=42]{Python}
house.add_component(
    technologies.HeatPump(
        thermal_power_limit=None
        cop_0_35=3.8,  # 1
    )
)
house.add_component(
    technologies.AirHeatExchanger(
        air_temperature=[3, 9],  # °C
    )
)
\end{minted}

In Figure~\ref{fig: example energy system graph},
the resulting energy system is displayed as a graph.
This graph was automatically generated by MTRESS building
upon using graphviz~\cite{Graphviz-overview}.

Once defined, the model is optimised using:
\begin{minted}[xleftmargin=\parindent,breaklines,linenos,firstnumber=53]{Python}
energy_system.optimise()
\end{minted}
Afterwards, results can be prompted
\begin{minted}[xleftmargin=\parindent,breaklines,linenos,firstnumber=53]{Python}
electricity_flows = energy_system.flows(Carrier, Electricity)
\end{minted}

The alternative is to define the energy system using a yaml file.
The following listing will result in exactly the same model as
the code explained above:
\begin{minted}[xleftmargin=\parindent,breaklines,linenos]{YAML}
general:
  timeindex:
    start: 2021-07-10 06:00:00
    end: 2021-07-10 08:00:00
    freq: 60T

locations:
  SFH:
    carriers:
      Heat:
        temperature_levels:
          - 20. # °C
          - 30. # °C
          - 55. # °C
        reference_temperature: 10 # °C
      Electricity:
        demand_rate: 0  # ct/kW
        working_price: 35  # ct/kWh

    demands:
      Electricity:
        name: "electricity demand"
        time_series: [7, 8.4]  # kW
      FixedTemperatureHeat:
        name: "space heating"
        time_series:
          - 13.37  # kW
          - 42  # kW
        flow_temperature: 30 # °C
        return_temperature: 20 # °C
      FixedTemperatureHeat:
        name: "hot water"
        time_series:
          - 0  # kW
          - 12  # kW
        flow_temperature: 55 # °C
        return_temperature: 10 # °C
    
    components:
      AirHeatExchanger:
        parameters:
          air_temperature: [3, 9]  # °C
      HeatPump:
        parameters:
          cop_0_35: 3.8
\end{minted}
Note that all comments are voluntary and are just meant
to help the reader.

For longer time series, it can be convenient not to inline every value.
This is supported by giving it in the form
\begin{minted}[xleftmargin=\parindent,breaklines,linenos,firstnumber=42]{YAML}
          air_temperature: file=weather.csv:temperature (°C)
\end{minted}
which denotes that the time series can be found
in the column ``temperature (\SI{}{\celsius})''
of the file ``weather.csv''.
Time series used in models can be either defined at points in time
or over time intervals.
For time-intervals, MTRESS uses left-indexed data,
meaning that the time stamp denotes the beginning
of the interval the value is valid.
Following the convention of solph this is being made
explicit by adding the time step limiting the last interval
with ``NaN'' (for ``not a number'') instead of data.
It might sound a bit pedantic at first,
but giving the length of the last interval is actually required because
there is no requirement that all time intervals have the same length.
Additionally, being explicit helps avoiding common off-by-one errors.

\section{Summary and Outlook}

We presented version~3.0 of MTRESS,
a tool that allows to create optimisation models for residential energy supply systems
in just a few lines of code.
These models can contain multiple temperature levels for the energy carrier heat
and multiple pressure levels for gases (like natural gas or hydrogen)
in a transparent and flexible way.
This way i.e. storage capacities can be modelled more realistically
compared to traditional linear models without levels.
There are two ways to create those models,
on one hand, there is a Python API for dynamic modelling,
on the other hand a YAML interface is provided.
For the future, it is planned to allow exporting models
as YAML files to facilitate reproducing results of energy system designs
that were created using a dynamic algorithm.

\section*{Acknowledgments}

We thank Uwe Krien for bringing up the topic of properly indexing data
and Lennart Schürmann for discussions about an improved storage formulation.
Elif Turhan and Herena Torio are acknowledged for proofreading the manuscript.

\bibliographystyle{unsrt}
\bibliography{osmses_2022-mtress_3.bib}

\end{document}